# ETUDE D'UN DISPOSITIF D'AIDE A L'INTENTION D'ELEVES EN DIFFICULTE DANS LA RESOLUTION D'UNE SITUATION-PROBLEME MATHEMATIQUE

Marie-Pier MORIN[*], Laurent THEIS[*], Teresa ASSUDE[**], Karine MILLON-FAURÉ[**], Jeannette TAMBONE[**], Jeanne KOUDOGBO[*], Valérie HAMEL[*], Patricia MARCHAND[*]

**Résumé** – Dans cette communication, nous analysons un dispositif d'aide mis en oeuvre par des enseignantes auprès d'élèves en difficulté suite à l'expérimentation de situations-problèmes mathématiques avec l'ensemble de la classe. À partir d'une première expérimentation (Theis et al., 2016), nous montrons d'abord les fonctions d'aide de ce dispositif, notamment à travers le triplet des genèses (Sensevy et al. 2000) : chronogenèse, mésogenèse et topogenèse. Ensuite nous analysons une deuxième mise en œuvre en détaillant, à partir d'épisodes concernant des élèves qui y ont participé, comment ces fonctions s'y retrouvent, ou non.

**Mots-clefs** : Situation-problème, élève en difficulté, dispositif d'aide, triplet des genèses, mathématiques

**Abstract** – In this paper, we analyze an aid session tested by an elementary school teachers. This aid session has been set up by a teacher for some students with difficulties after the work in the whole class. We first show how this aid session can help the pupils, especially through its chrono-topo-meso-genetic functions. Then we analyze a second situation in which the aid session is implemented, looking for clues of these functions in some episodes involving students who participated in the aid session.

Notre communication se situe dans le cadre du groupe 11 et concerne les dispositifs pour la scolarisation des publics spécifiques, plus précisément les dispositifs mis en oeuvre par les enseignants pour aider les élèves qu'ils jugent en difficulté. Dans cette perspective, nous étudierons un dispositif d'aide post intervention en classe mis en place par des enseignantes en portant notre attention sur les fonctions topogénétiques, chronogénétiques et mésogénétiques de ce dispositif. À cette fin, nous présenterons notre cadre théorique et problématique à partir d'une première expérimentation de ce dispositif (Theis et al., 2016) et analyserons par la suite une nouvelle expérimentation de ce dispositif auprès de quatre élèves identifiés en difficulté par leur enseignante. Enfin, nous tenterons de voir si les effets rencontrés dans le premier dispositif analysé se retrouvent dans le dispositif que nous présentons.

## I. CADRE THEORIQUE ET PROBLEMATIQUE

Lors de l'implantation du renouveau pédagogique, dont les orientations sont décrites dans le Programme de formation de l'école québécoise (Gouvernement du Québec, 2001), une place importante a été réservée à la résolution de situations-problèmes. D'une part, la compétence à résoudre une situation-problème est centrale dans le programme de mathématiques, tant au primaire qu'au secondaire, et d'autre part la compétence consistant à résoudre un problème est transversale au programme. Aussi, le Ministère de l'Éducation (MEQ) considère la situation-problème comme « un outil intellectuel puissant au service du raisonnement et de l'intuition créatrice » (p. 126). Par conséquent, la résolution de situations-problèmes doit être abordée comme objet d'apprentissage et comme outil dans l'apprentissage des mathématiques. Pour les élèves qui présentent des difficultés, la résolution de situations-problèmes constitue ainsi un défi important pouvant même conduire à l'échec. C'est dans ce cadre que nous avons mené deux recherches collaboratives avec des enseignantes du primaire, lesquelles portaient sur les conditions favorables à l'engagement dans une situation-problème et l'apprentissage des concepts visés. La première a été réalisée auprès de huit enseignantes d'une école primaire de Sherbrooke et la seconde, auprès de huit enseignantes de plusieurs écoles primaires situées en périphérie de Sherbrooke. Dans le cadre de ces recherches collaboratives, nous avons accompagné les enseignantes dans l'élaboration, l'expérimentation et l'analyse de situations-problèmes mathématiques en classe. Ces situations-problèmes, construites autour des critères d'une situation-problème d'Astolfi (1993) étaient pensées de manière à permettre aux élèves de faire de nouveaux apprentissages et de surmonter des obstacles de nature épistémologique.

---

[*] Université de Sherbrooke – Québec – Canada – Marie-Pier.Morin@USherbrooke.ca; Laurent.Theis@USherbrooke.ca; Jeanne.Koudogbo@USherbrooke.ca; Valerie.Hamel@USherbrooke.ca; Patricia.Marchand@USherbrooke.ca.

[**] ADEF – Université d'Aix-Marseille – France – teresa.dos-reis-assude@univ-amu.fr; karine.millon-faure@univ-amu.fr; jane.tambone@wanadoo.fr.

Dans la première phase du projet, Sylvie, une enseignante du $2^e$ cycle du primaire[1], a mis sur pied un dispositif particulier, qui avait comme but de rencontrer les élèves qu'elle jugeait en difficulté avant la résolution de la situation-problème en classe. Ce dispositif visait à permettre à ces élèves de faire connaissance avec la situation-problème travaillée en classe sans toutefois déjà y travailler sur les objets de savoir visés. L'analyse de l'expérimentation de cette mesure (Theis et al., 2014), a fait ressortir trois fonctions potentielles de ce dispositif d'aide pour les élèves en difficulté. Par la suite, de nouvelles versions du dispositif, toutes différentes de la première, ont été expérimentées et analysées (Assude et al., 2015, 2016; Millon-Fauré et al., sous presse; Theis et al., 2016).

Nous avons modélisé les fonctions de notre dispositif d'aide autour du triplet des genèses décrit par Sensevy, Mercier et Schubauer-Leoni (2000) pour qui,

> « au sein du système didactique, le professeur doit agir (définir, réguler, dévoluer, instituer) pour : produire les lieux du professeur et de l'élève (effet de topogénèse); produire les temps de l'enseignement et de l'apprentissage (effet de chronogénèse); produire les objets des milieux des situations et l'organisation des rapports à ces objets (effet de mésogénèse) » (p. 267).

Parmi les trois fonctions identifiées, la fonction mésogénétique permet la rencontre de l'élève en difficulté avec différents éléments du milieu, avant la rencontre en classe de la situation-problème. La fonction chronogénétique permet à l'élève d'en savoir davantage sur la tâche à réaliser, avant les autres élèves. L'élève dispose également de plus de temps pour entrer dans la tâche. Toutefois, nos analyses ont révélé que le temps didactique n'avançait pas dans le dispositif, puisqu'aucune institutionnalisation ne s'y produit. Cependant, l'élève en apprend davantage sur les tâches à réaliser et sur les techniques qui y seront pertinentes, ainsi que sur les discours les justifiant. Nous avons qualifié cette progression d'avancée du temps praxéologique (Assude et al.; 2016) Finalement, la fonction topogénétique du dispositif « avant » permet à l'élève en difficulté d'assumer pleinement son rôle d'élève lors du travail en classe.

---

[1] Élèves âgés de 8 à 10 ans.

Toujours au cours de la première phase du projet, Manon, une enseignante du 3$^e$ cycle du primaire[2], a expérimenté ce dispositif en y intégrant un dispositif d'aide supplémentaire, postérieurement à la situation vécue en classe. Ce dispositif a également fait l'objet d'analyses (Theis et al., 2016), de façon à voir si les fonctions identifiées lors de l'expérimentation de la mesure d'aide avant la situation-problème se retrouvent également lors de l'expérimentation de la mesure d'aide suite à la situation-problème.

Dans le cadre de cette communication, à partir de l'analyse du dispositif post de Manon, nous analyserons la mise en œuvre d'un deuxième dispositif post réalisé cette fois-ci avec la classe d'Isabelle, et verrons en quoi les fonctions qui se dégagent du dispositif mis en place par Manon, se retrouvent, ou non, dans celui expérimenté par Isabelle.

Pour faire l'analyse de ce dispositif, à l'instar de Chevallard (1995), nous modéliserons ce regroupement d'élèves considérés en difficulté comme un système didactique auxiliaire (SDA) à la classe. Ce dernier considère l'espace d'étude comme un ensemble de systèmes didactiques principaux (SDP), la classe par exemple, et des SDA, qui aident à l'étude. Les SDA étant entre autres liés aux mêmes enjeux de savoirs que les SDP, ils dépendent des SDP. À l'image des SDA décrits ailleurs (Tambone, 2014), le SDA que nous étudierons a été réalisé en aval de la situation-problème vécue en classe (SDA post). Ainsi, le SDA post fait suite au SDP, qui lui-même se situe après le SDA pré. Ces dispositifs d'aide (pré et post) peuvent être illustrés par le schéma qui suit :

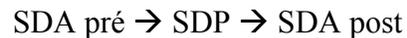

*Figure 1 – Séquence des dispositifs d'aide*

Dans la prochaine section, à partir de l'analyse de la mise en place du SDA post de Manon, nous pourrons apprécier la façon dont ces fonctions s'actualisent.

*1. Analyse du SDA post de Manon*

Le SDA post de Manon, pour qui le SDP (Theis et al, 2016) portait sur le calcul de l'aire du triangle, comportait les sept étapes suivantes :

1. Justification du dispositif auprès des élèves.
2. Vérification, de manière individuelle, auprès de chaque élève s'ils seraient maintenant en mesure de calculer l'aire du triangle.
3. Travail sur l'aire du triangle isocèle utilisé dans le SDP à travers la technique « inscription dans un rectangle, puis diviser l'aire du rectangle par 2 ».
4. Travail sur un cas particulier, en utilisant la technique « trouver différents endroits où on peut placer la hauteur ».
5. Nouvelle vérification par Manon si les élèves seraient maintenant en mesure de déterminer seuls l'aire d'un triangle.
6. Retour sur le calcul de l'aire du rectangle à travers des exemples utilisés dans le SDP et lien avec l'aire du triangle.
7. Manon questionne les élèves sur la façon dont ils ont perçu l'utilité du SDA pré.

---

[2] Élèves âgés de 10 à 12 ans.

L'analyse de la mise en place de ces étapes à l'aide du triplet des genèses a révélé des fonctions potentielles du dispositif. Premièrement, la fonction mésogénétique s'est manifestée à travers des objets qui sont apparus dans le SDA post, mais qui étaient restés implicites avec l'ensemble du groupe dans le cadre du SDP. Par exemple, lors de la 3$^e$ étape, Manon est revenue avec une élève sur l'utilité de la formule de l'aire du triangle, élément qui n'avait pas été abordé de façon explicite avec le groupe. Aussi, l'analyse a montré que le mode de fonctionnement en groupe restreint dans le cadre du SDA a permis de clarifier des objets encore problématiques pour certains élèves. Enfin, toujours au niveau mésogénétique, la séance post a permis à Manon de revenir sur des objets plus anciens, comme la compréhension de la formule de l'aire du rectangle, nécessaire à la compréhension de la formule de l'aire du triangle.

Deuxièmement, au plan topogénétique, l'intervention en petit groupe dans le cadre du SDA post a permis à l'enseignante d'intervenir avec les élèves de façon individuelle, comme le montre l'étape 2, où Manon prend du temps avec chacun des élèves pour vérifier leur compréhension du calcul de l'aire d'un triangle. Également, nous avons pu observer que le fonctionnement en groupe restreint a permis aux élèves de poser des questions d'éclaircissement à Manon. C'était entre autres le cas de Roselyne qui se questionnait sur le besoin de diviser l'aire par 2 lorsque la hauteur est inscrite dans le triangle.

Enfin, au plan chronogénétique, l'analyse de la séance a montré que les élèves ciblés disposent de plus de temps pour l'institutionnalisation. Ce temps supplémentaire a permis à Manon de revenir sur des éléments restés implicites pour les élèves. Également, nous avons noté que, étant donné que le SDA post de Manon s'est déroulé une semaine après l'institutionnalisation en classe dans le SDP, cette séance est apparue comme une deuxième institutionnalisation – une institutionnalisation bis –, où Manon a convoqué à nouveau la technique institutionnalisée, c'est-à-dire la formule de l'aire du triangle.

## II. DEROULEMENT DU SDA POST D'ISABELLE

Le dispositif principal de la situation-problème proposée par Isabelle porte sur les fractions et leur écriture. Les élèves, qui sont au 2$^e$ cycle du primaire[3], doivent partager des quantités et exprimer à l'écrit leurs réponses avec des fractions et des nombres fractionnaires[4]. Comme première tâche, à partir de cercles en carton et de ciseaux, les élèves doivent partager également 5 fromages entre 2 élèves. Comme deuxième tâche, toujours à partir de matériel concret, les élèves doivent partager également 7 tartes entre 4 enfants. La troisième tâche consiste à partager 4 carrés également entre 8 enfants. Contrairement aux deux premières tâches, celle-ci doit être réalisée entièrement sur une tablette numérique. Comme dernière tâche, à partir d'une feuille sur laquelle sont dessinés 6 rectangles, les élèves doivent partager 6 tablettes de chocolat entre 8 personnes. Pour toutes les tâches, la codification doit être faite sur une tablette numérique, de façon à pouvoir projeter les réponses au tableau blanc par la suite.

Lors du SDA post, l'enseignante propose des tâches qui, pense-t-elle, lui permettent de revisiter ce qui a été travaillé en classe. Trois tâches sont prévues : représenter d'une manière analogique le nombre fractionnaire 2 ½ ; partager 2 barres de chocolat entre 3 personnes ; partager 3 barres de chocolat entre 5 personnes.

---
[3] 8 à 10 ans
[4] Cette situation est inspirée des recommandations de Van de Walle et Lovin (2008).

Pour être en mesure de mieux comprendre la séance, nous dressons la description suivante du SDA post :

1. Appel à la mémoire du groupe concernant l'activité vécue dans le cadre du SDP.
   « Ok, donc on aujourd'hui on se revoit. On a fait notre activité vendredi et aujourd'hui on est lundi. Puis moi j'aimerais savoir, comme ça, qu'est-ce que tu as retenu de cette activité-là ? »

2. Vérification de la compréhension du concept de la fraction auprès des 4 élèves ciblés.

3. Présentation d'une première tâche : représenter le nombre fractionnaire 2 ½.
   « On va y aller un petit peu à l'envers de ce qu'on a fait [Isabelle écrit sur des feuilles le nombre fractionnaire 2 ½] Es-tu capable de me l'illustrer Jules ? [Isabelle distribue au fur et à mesure les feuilles] ».
   Chacun des élèves essaie de représenter le nombre. Lors du retour collectif, Isabelle se rend compte que les élèves ont de la difficulté à représenter la partie entière.

4. Présentation de la deuxième tâche : partager 2 entiers (rectangles dessinés) en 3.
   « On va refaire un peu comme l'exercice qu'on a fait en classe [Isabelle distribue les feuilles], ok. Puis là, il va falloir que tu m'écrives qu'est-ce que ça fait. Donc, j'ai 2 barres de chocolat. Et moi, je veux les partager en 3 ».
   Chacun des élèves essaie de représenter le partage et de codifier la réponse. Tous réutilisent la technique utilisée dans le SDA pré, soit de partager chacun des entiers en 3 et de numéroter les 6 morceaux (1, 2, 3, 1, 2, 3). Alors que les élèves éprouvent des difficultés à codifier le partage, Isabelle les guide vers l'écriture fractionnaire.

5. Présentation de la troisième tâche : partager 3 entiers (rectangles dessinés) en 5.
   « J'ai 3 barres de chocolat, puis mes 3 barres de chocolat, je veux les partager en 5 personnes ».
   Chacun des élèves essaie de représenter le partage et de codifier la réponse. Isabelle réalise un retour individuel pour deux élèves et un retour collectif pour les deux autres qui ont plus de difficulté. La codification reste difficile pour les quatre élèves, particulièrement pour ceux qui ont commencé par diviser chaque rectangle en 2 pour faire ½ par personne, et qui ont ensuite divisé à nouveau la demie restante, pour faire 1/10 par personne.

III. ANALYSE DU SDA POST D'ISABELLE SELON LES FONCTIONS DECRITES

   *1. Fonction mésogénétique*

Dans le cadre du SDA post d'Isabelle, nous avons pu retrouver des effets qui se sont manifestés dans le SDA post de Manon. En effet, au plan mésogénétique, nous avons noté que la séance post a permis aux élèves en difficulté de revenir sur le milieu et de clarifier certains objets encore sensibles comme le partage en parts égales et la codification de la fraction. Par contre, la représentation du nombre fractionnaire, qui n'a pas été approfondie dans le SDP, n'est pas allée en ce sens. En effet, dans le cadre du SDP, même si deux réponses sur quatre impliquaient un nombre fractionnaire écrit sous la forme a b/c (2 ½ et 1 ¾), l'accent a plutôt été mis sur la façon de codifier les parties fractionnaires.

Ainsi, pour débuter le SDA post, l'enseignante a demandé aux élèves d'illustrer 2 ½. L'entrevue qui a précédé la séance post permet de penser qu'Isabelle n'avait pas prévu de demander de faire illustrer un nombre fractionnaire par les élèves et que cette décision a été prise dans l'action, suite à un commentaire d'Olivier au tout début du SDA post :

> « Isabelle : Fiou ! ça été dur. Donc ça une fraction et ça c'est une fraction [Isabelle pointe les fractions ¾ et 1/3]. On sait que dans la fraction on a toujours deux parties. Les chiffres qui sont en haut et les chiffres qui sont en bas.
> Olivier : Et à côté ?
> Isabelle : Hein ?
> Olivier : Et à côté ?
> Isabelle : Des fois il y en a à côté. Qu'est-ce que tu veux dire ? C'est comme si je mettais un 2 là. Comme ça? [Isabelle écrit 2 sur la feuille] Puis la fraction tu l'écrirais comment ? À côté. Si je te dis une demie, écris-moi là à côté [Isabelle remet la feuille à Olivier].
> [L'élève a un peu de difficulté à écrire ½, mais après hésitations, il réussit]
> Isabelle : Ok, si je voulais illustrer ça : 2 et ½. Qu'est-ce que je ferais ? Est-ce que vous êtes capable de me le faire tout le monde ? »

Isabelle a toutefois introduit ce nouveau type de tâche comme étant le contraire de ce qui avait été fait auparavant :

> « On va y aller un petit peu à l'envers de ce qu'on a fait [Isabelle écrit sur des feuilles le nombre fractionnaire 2 ½]. »

Il ne s'agit pas du même type de tâche qu'elle propose aux élèves. L'enseignante fait appel à la mémoire didactique des élèves en convoquant un objet qu'elle croyait ancien. Par contre, s'il s'agissait d'un objet ancien s'inscrivant dans la continuité de la tâche pour elle, compte tenu qu'il ne s'agit pas du même type de tâche, cet objet est devenu nouveau pour les élèves.

Aussi, concernant les nombres choisis pour le SDA post, puisqu'ils impliquaient un partage en un nombre impair (3 et 5), contrairement aux partages faits lors du SDP (2, 4 et 8), nous considérons que le niveau de difficulté n'était pas le même, donc que le milieu était un peu différent. Par contre, ces nombres ont été choisis de façon éclairée par Isabelle.

> « Et je vais leur redonner un petit problème avec un rectangle qu'ils vont devoir fractionner en n'utilisant pas le 2. Justement pour voir comment ils vont se débrouiller avec le partage à ce moment-là. Où ils ne pourront pas y aller instinctivement : bien là je vais le découper en 2, puis là ça ne marche pas, je vais y aller en 4. Donc, je veux vérifier ça avec eux. Parce que ça, parce qu'ils n'ont pas été capables d'aller au-delà de ça. »

Les élèves ont donc été placés dans un milieu un peu différent de ce dont ils ont connu en SDP, mais cette modification était calculée par l'enseignante, qui souhaitait provoquer une technique de fractionnement différente lors du partage de l'entier.

2. *Fonction topogénétique*

Au niveau topogénétique, l'analyse permet de voir que les quatre élèves sont très engagés dans le dispositif et prennent réellement leur place d'élèves. Pour chacune des tâches, l'enseignante est revenue avec chacun des élèves, en considérant les différentes techniques présentées et en allouant plus de temps à ceux qui en avaient plus besoin. On peut donc dire qu'au niveau conceptuel, elle laisse de l'espace et de l'ouverture aux élèves. Cela dit, lors de la deuxième tâche, comme on peut le voir dans l'extrait suivant, elle était somme toute assez près du guidage pour aider les élèves à codifier leur fraction.

Aussi, le groupe restreint a permis aux élèves de poser des questions et de faire des affirmations qu'ils n'auraient peut-être pas faites en grand groupe. Prenons comme exemple le cas d'Olivier qui, suite à une question d'Isabelle, a affirmé qu'il peut y avoir un nombre à côté de la fraction. Comme nous venons de le voir, cette affirmation a fait en sorte qu'Isabelle a fait le choix d'introduire un nouvel objet dans le milieu.

### 3. Fonction chronogénétique

Au plan chronogénétique, la séance post a permis aux élèves d'avoir plus de temps que les autres élèves pour approfondir leur compréhension d'objets implicites ou encore, d'objets qui étaient restés sensibles, comme la codification des fractions et le partage en parts égales. Par exemple, dans l'extrait suivant, on voit que l'enseignante prend du temps pour distinguer la fonction du numérateur et du dénominateur, sans toutefois nommer aucun des deux, tout comme ça été le cas lors du SDP.

> « Ok. Comment tu vas écrire ta réponse pour me dire combien ça faire de morceaux que chacun reçoit ? »
> « 2 »
> « 2 morceaux. Donc ça, s'il reçoit 2 morceaux ça va aller où dans ta fraction ? En haut ou en bas ? »
> « En haut. »
> « En haut. Ok, écris-le. [Olivier prend son crayon] Écris la partie en haut puis après on va regarder comment faire pour la partie bas »
> « OK »
> « Qu'est-ce qu'on va écrire en bas ? C'est quoi mon tout ? »

Par contre, contrairement au SDA post de Manon, nous ne pouvons pas affirmer qu'il y a eu institutionnalisation forte puisque l'enseignante s'est à peine détachée du contexte pour prendre une distance. Nous pouvons constater la même chose à l'aide de l'exemple suivant, qui concerne le partage en parts égales :

> « Penses-tu qu'il va y avoir un ami jaloux si tu lui donnes un gros morceau comme ça et puis les deux autres tu ne leur donnes juste une demie, un petit morceau? ».

Ainsi, Isabelle reste toujours très près du contexte des problèmes présentés, ce qui ne nous permet pas de dire qu'il y a eu une reprise de l'institutionnalisation, contrairement au SDA post de Manon.

## IV. CONCLUSION

En conclusion, malgré les quelques différences observées, nous avons pu faire ressortir des éléments communs aux deux SDA post. Au plan mésogénétique on peut dire que dans les deux cas il y a eu retour sur des objets encore sensibles. Par contre, la séance post de Manon lui avait permis de revenir sur des objets restés implicites ainsi que sur des objets plus anciens, ce qui ne semble pas avoir été le cas pour Isabelle. De même, cette enseignante a introduit un nouveau type de tâche dans le milieu, ce qui s'est avéré problématique. Au plan topogénétique, dans les deux cas, le travail en petit groupe que permet le SDA a permis aux enseignantes de donner plus de temps à leurs élèves et aussi, a fait en sorte que les élèves pouvaient poser plus de questions à leur enseignante. De plus, dans les deux SDA post, on peut dire que les élèves étaient très engagés. Les différences sont plus importantes au plan chronogénétique. En effet, dans les deux cas il y a eu plus de temps pour les apprentissages, mais dans le cas de Manon, il s'agissait d'une institutionnalisation bis, ce qui n'a pas été le cas pour Isabelle.

À la lumière de ces analyses nous pouvons affirmer que nous retrouvons des effets potentiels des trois fonctions évoquées plus haut dans l'expérimentation des deux dispositifs post. Évidemment, cela prendra d'autres analyses pour pouvoir tirer des conclusions en ce sens. Ces analyses, qui sont en cours de réalisation, seront présentées lors de notre communication.